\documentclass[11pt]{article}
\usepackage[dvips]{graphicx}
\usepackage{color}
\usepackage{amsmath}
\usepackage{amsfonts}
\usepackage{amssymb}
\usepackage{amsthm}
\usepackage{newlfont}
\usepackage{epsfig}
\usepackage{multirow}
\usepackage{setspace}
\usepackage{hyperref}
\usepackage{cite}

\begin{document}

\newtheorem{assumption}{Assumption}[section]
\newtheorem{definition}{Definition}[section]
\newtheorem{lemma}{Lemma}[section]
\newtheorem{proposition}{Proposition}[section]
\newtheorem{theorem}{Theorem}[section]
\newtheorem{corollary}{Corollary}[section]
\newtheorem{remark}{Remark}[section]

\title{Computing Weakly Reversible Linearly Conjugate Chemical Reaction Networks with Minimal Deficiency}
\author{Matthew D. Johnston$^a$, David Siegel$^a$ and G\'{a}bor Szederk\'{e}nyi$^{b,c}$ \bigskip \\
${}^a$ Department of Applied Mathematics,\\
University of Waterloo,\\
Waterloo, Ontario, Canada N2L 3G1\\
${}^b$ Systems and Control Laboratory,\\
Computer and Automation Research Institute,\\
Hungarian Academy of Sciences\\
H-1518, P.O. Box 63, Budapest, Hungary\\
${}^c$ Faculty of Information Technology,\\
P\'{e}ter P\'{a}zm\'{a}ny Catholic University,\\
H-1083, P.O. Box 278, Budapest, Hungary}
\date{}
\maketitle

\tableofcontents

\bigskip

\begin{abstract}

Mass-action kinetics is frequently used in systems biology to model the behaviour of interacting chemical species. Many important dynamical properties are known to hold for such systems if they are weakly reversible and have a low deficiency. In particular, the Deficiency Zero and Deficiency One Theorems guarantee strong regularity with regards to the number and stability of positive equilibrium states. It is also known that chemical reaction networks with disparate reaction structure can exhibit the same qualitative dynamics. The theory of linear conjugacy encapsulates the cases where this relationship is captured by a linear transformation. In this paper, we propose a mixed-integer linear programming algorithm capable of determining weakly reversible reaction networks with a minimal deficiency which are linearly conjugate to a given reaction network.

\end{abstract}


\bigskip

\section{Introduction}

A chemical reaction network is given by sets of chemical reactants reacting at prescribed rates to form sets of chemical products. Under suitable assumptions, such as mass-action kinetics and spatial homogeneity, the time evolution of the concentrations of the chemical species can be modeled by a set of autonomous polynomial ordinary differential equations. Such mass-action systems are frequently used to model systems in systems biology and other areas of computational biology \cite{Sh-F,T-G,T-G2}.

The systematic study of chemical reaction networks was initiated in 1972 in the papers \cite{F1,H,H-J1}. In \cite{H-J1}, the authors presented a condition on the positive equilibrium concentrations, called \emph{complex balancing}, which is sufficient to guarantee that within each linear invariant space of the network there exists a unique equilibrium concentration, and that this concentration is locally asymptotically stable relative to its invariant space. In \cite{F1} and \cite{H}, the authors related the capacity of a network to exhibit complex balancing to a nonnegative network parameter called the \emph{deficiency}. In particular, they showed that a network is complex balanced for all sets of rate constant values if and only if it is weakly reversible and has a deficiency of zero. This \emph{deficiency}-oriented approach to analysing chemical reaction networks has since been applied to a wide variety of biochemical systems, including enzymatic models, signal transduction, and phosphorylation networks \cite{S-M,S,M-G,D-M,M-D-S-C}.

It is also known that many qualitative properties of the dynamics of reaction networks can be shared by networks with disparate network structure. A thorough study of \emph{dynamical equivalence}, the property that two reaction networks give rise to identical mass-action dynamics, was conducted in \cite{C-P}. In \cite{J-S2}, dynamical equivalence was extended to \emph{linear conjugacy}, whereby the trajectories of two mass-action systems could be related by a non-trivial linear transformation. The problem of determining dynamically equivalent networks with the greatest and fewest number of reactions was placed in a mixed-integer linear programming (MILP) context in \cite{Sz2}. This methodology has been extended to linear conjugacy in the papers \cite{J-S4} and \cite{J-S5}.

Linear conjugacy effectively creates classes of mechanisms which, despite disparate network structure and properties, have equivalent dynamics. A long-standing problem in chemical reaction network theory, first stated in \cite{H-T}, is to ``... look for a mechanism in a class of mechanisms with a given - chemically relevant - property. Such a property may be conservativity, (weak) reversibility, zero deficiency or just structural stability as well.'' In this paper, we use linear conjugacy and MILP methods to answer this challenge for the class of weakly reversible mechanisms where the structural property of interest is the deficiency.

In general, we will be interested in mechanisms within the class of linearly conjugate networks with \emph{minimal} deficiency. This is due to the observation that the behaviour of a chemical reaction network's mass-action system is generally more regular (e.g. fewer steady states, lower capacity for oscillations, etc.) for networks with a lower deficiency than a higher deficiency \cite{F1,H,F2,Fe4}. Consequently, within the class of mechanisms which are linearly conjugate, if there are mechanisms with differing deficiencies, the network with the lowest deficiency is generally preferred. We therefore modify the existing MILP framework to compute weakly reversible linearly conjugate networks which have a minimal deficiency. The methodology is based on known properties of the kernel of the kinetics matrix of a weakly reversible chemical reaction network.

\section{Background}

In this section we present terminology and notation relevant to the study of chemical reaction networks. We introduce the notion of the deficiency of a network and a few classical results which relate the deficiency of a network to the dynamics of the network's corresponding mass-action system. We also introduce the notion of two chemical reaction networks being linearly conjugate.

\subsection{Chemical Reaction Networks}

The chemical \emph{species} or \emph{reactants} of a network will be given by the set $\mathcal{S} = \left\{ X_1, X_2, \ldots, X_n \right\}$. The combined elements on the left-hand and right-hand side of a reaction are given by linear combinations of these species. These combined terms are called \emph{complexes} and will be denoted by the set $\mathcal{C} = \left\{ C_1, C_2, \ldots, C_m \right\}$ where
\[C_i = \sum_{j=1}^n \alpha_{ij} X_j, \; \; \; i=1, \ldots, m\]
and the $\alpha_{ij}$ are nonnegative integers called the \emph{stoichiometric coefficients}. The complex with all stoichiometric coefficients equal to zero will be called the \emph{null complex} and denoted by $C_i = \emptyset$. We define the reaction set to be $\mathcal{R} = \left\{ (C_i,C_j) \; | \; C_i \mbox{ reacts to form } C_j \right\}$. The property $(C_i,C_j) \in \mathcal{R}$ will commonly be denoted $C_i \to C_j$. To each $(C_i,C_j) \in \mathcal{R}$ we will associate a positive \emph{rate constant} $k(i,j) > 0$ and to each $(C_i,C_j) \not\in \mathcal{R}$ we will set $k(i,j) = 0$. The triplet $\mathcal{N} = (\mathcal{S}, \mathcal{C}, \mathcal{R})$ will be called the \emph{chemical reaction network}.

The above formulation naturally gives rise to a directed graph $G(V,E)$ where the set of vertices is given by $V = \mathcal{C}$, the set of directed edges is given by $E = \mathcal{R}$, and the rate constants $k(i,j)$ correspond to the weights of the edges from $C_i$ to $C_j$. In the literature this has been termed the \emph{reaction graph} of the network \cite{H-J1}. Since complexes may be involved in more than one reaction, as a product or a reactant, there is further graph theory we may consider. A \emph{linkage class} is a maximally connected set of complexes, that is to say, two complexes are in the same linkage class if and only if there is a sequence of reactions in the reaction graph (of either direction) which connects them. We will denote by $\ell$ the number of linkage classes in a network. A reaction network is called \emph{weakly reversible} if $C_i \to C_j$ for any $C_i, C_j \in \mathcal{C}$ implies there is some sequence of complexes such that $C_j = C_{\mu(1)} \to C_{\mu(2)} \to \cdots \to C_{\mu(l-1)} \to C_{\mu(l)} = C_i$.


A directed graph is called \emph{strongly connected} if there exists a directed path from each vertex to every other vertex. A \emph{strongly connected component} of a directed graph is a maximal set of vertices for which paths exist from each vertex in the set to every other vertex in the set. A strongly connected component is called \emph{terminal} if there is no reaction leading from a vertex in the strongly connected component to a vertex not in the component. For a weakly reversible network, all strongly connected components are terminal, and they correspond to the linkage classes of the reaction graph.

Assuming mass-action kinetics, the dynamics of the specie concentrations over time is governed by the set of differential equations
\begin{equation}
\label{de}
\frac{d\mathbf{x}}{dt} = Y \cdot A_k \cdot \Psi(\mathbf{x})
\end{equation}
where $\mathbf{x} = [ x_1 \; x_2 \; \cdots \; x_n ]^T$ is the vector of reactant concentrations. The \emph{stoichiometric matrix} $Y$ contains entries $[Y]_{ij} = \alpha_{ji}$ and the \emph{Kirchhoff} or \emph{kinetics} matrix $A_k$ is given by
\begin{equation}
\label{kinetics}
[A_k]_{ij} = \left\{ \begin{array}{cll} -\sum_{l=1,l \not= i}^m k(i,l), & \mbox{  if  } & i = j \\ k(j,i) & \mbox{  if  } & i \not= j \end{array} \right.
\end{equation}
for $i,j=1, \ldots, m$. When we speak of the \emph{structure} of a kinetics matrix, we are referring to the distribution of positive and zero entries, which determines the network structure of the corresponding reaction graph. Finally, the \emph{mass-action vector} $\Psi(\mathbf{x})$ is given by
\begin{equation}
\label{psi}
\Psi_j(\mathbf{x}) = \prod_{i=1}^n x_i^{[Y]_{ij}}, \; \; \; j=1, \ldots, m.
\end{equation}

The behaviour of solutions of (\ref{de}) can be further understood by consideration of the \emph{reaction vectors}
\[v_{ij} = \left\{ \begin{array}{ll} [Y]_{\cdot,j} - [Y]_{\cdot,i} \; \; \; \; \; & \mbox{for } (C_i, C_j) \in \mathcal{R} \\ 0 & \mbox{otherwise} \end{array} \right.\]
where $[Y]_{\cdot,i}$ denotes the $i$th column of $Y$. The span of the vectors is called the \emph{stoichiometric subspace} and is denoted by $S=$ span$\left\{ v_{ij} \; | \; (C_i,C_j) \in \mathcal{R} \right\}$. We will denote the dimension of $S$ by $s=$ dim$(S)$.

It is clear that the right-hand side of (\ref{de}) is contained in $S$ so that the vector field always directs trajectories within an affine translation of $S$. It can be shown that trajectories are restricted to the \emph{stoichiometric compatibility classes} $(\mathbf{x}_0 + S) \cap \mathbb{R}^n_{>0}$ \cite{H-J1,V-H}.

\subsection{Complex Balanced Networks and Deficiency}
\label{deficiencysection}

The structure of the reaction graph of a chemical reaction network often plays an important role in determining the dynamical behaviour of trajectories of (\ref{de}).

A particularly important structural parameter of a chemical reaction network is the \emph{deficiency}, which was introduced in \cite{H,F1} and has been studied extensively since \cite{Fe2,Fe3,Fe4,S-C}.
\begin{definition}
\label{deficiency}
The \textbf{deficiency} of a chemical reaction network is given by
\[\delta = m - \ell - s\]
where $m$ is the number of stoichiometrically distinct complexes, $\ell$ is the number of linkage classes, and $s$ is the dimension of the stoichiometric subspace.
\end{definition}
\noindent It is easy to show that the deficiency of a network may only take on nonnegative integer values.

The deficiency of a network is often related to properties of the equilibrium set of the mass-action system (\ref{de}). It is strongly related to the capacity for a network to admit \emph{complex balanced equilibrium concentrations}.

\begin{definition}
An equilibrium concentration $\mathbf{x}^* \in \mathbb{R}_{>0}^n$ is called a \textbf{complex balanced equilibrium concentration} if
\[A_k \cdot \Psi(\mathbf{x}^*) = \mathbf{0}.\]
A chemical reaction network is called \textbf{complex balanced} if every equilibrium concentration is a complex balanced equilibrium concentration.
\end{definition}
\noindent 
It is known that if a network is complex balanced at one equilibrium concentration then it is complex balanced at all of them (Theorem 6A, \cite{H-J1}). Consequently, a network for which any equilibrium concentration is complex balanced is a complex balanced network.

The following results relate the deficiency of a network to its capacity for complex balancing (see Theorem 4A of \cite{H} and Theorem 4.1 of \cite{F1}).
\begin{theorem}[Deficiency Zero Theorem]
\label{deficiencyzerotheorem}
A chemical reaction network is complex balanced for all rate constant values if and only if it is weakly reversible and has a deficiency of zero.
\end{theorem}

\begin{theorem}
\label{conditionstheorem}
If a mass-action system is weakly reversible, then the deficiency corresponds to the number of algebraically independent conditions on the rate constants which need to be satisfied in order for the system to be complex balanced.
\end{theorem}
\noindent Theorem \ref{deficiencyzerotheorem} is a special case of Theorem \ref{conditionstheorem}, taking the deficiency to be zero. We state is separately, however, due to its historical importance.

Strong dynamical properties are known to follow from complex balancing; in particular, the following was proved in \cite{H-J1}.
\begin{theorem}[Theorem 6A and Lemma 4C, \cite{H-J1}]
\label{cbtheorem}
If a chemical reaction network is complex balanced then there exists within each positive stoichiometric compatibility class exactly one equilibrium concentration, and that equilibrium concentration is locally asymptotically stable relative to its compatibility class.
\end{theorem}

The surprising implication of Theorem \ref{deficiencyzerotheorem} and Theorem \ref{cbtheorem} is that we can know very strong properties about the equilibrium set of the mass-action system (\ref{de}) based on structural properties of the \emph{reaction graph alone}. Most strikingly, the results hold regardless of the values of the rate constants. This is particularly important for biological examples, where the relevant rate constants are frequently unknown or unknowable. In cases where the deficiency in nonzero, supplemental conditions on the rate constants are required to bridge the gap from weak reversibility to complex balancing. Theorem \ref{conditionstheorem} tells us that we are more likely to be interested in networks with a \emph{lower} deficiency, since networks with a lower deficiency are, in some senses, \emph{closer} to complex balancing (and the regular dynamics guaranteed by Theorem \ref{cbtheorem}) than networks with a higher deficiency.

A result which further relates deficiency to properties of the equilibrium set of a chemical reaction network is the following, which can be found in \cite{F2}. It is worth noting, once again, that a lower deficiency generally guarantees more predictable behaviour than a higher deficiency. (Further work investigating the capacity of a deficiency one chemical reaction network to exhibit multistability was presented in \cite{Fe4}.)

\begin{theorem}[Deficiency One Theorem]
\label{deficiencyonetheorem}
Consider a chemical reaction network with deficiency $\delta$ and linkage classes $\mathcal{L}_1, \ldots, \mathcal{L}_\ell$. Let $\delta_i$, $i=1, \ldots, \ell$, be the deficiency of the linkage class $\mathcal{L}_i$ considered as its own subnetwork. Suppose that
\begin{enumerate}
\item
$\delta_i \leq 1, \mbox{ for } i=1, \ldots, \ell;$
\item
$\sum_{i=1}^\ell \delta_i = \delta; \mbox{ and}$
\item
each linkage class contains only one terminal strongly connected component.
\end{enumerate}
Then, if the network admits a positive equilibrium concentration, there is exactly one equilibrium concentration in each positive stoichiometric compatibility class. Furthermore, if the network is weakly reversible, then for every set of rate constants the network admits a positive equilibrium concentration.
\end{theorem}

\subsection{Linearly Conjugate Networks}

Under the assumption of mass-action kinetics, it is possible for the trajectories of two reaction networks $\mathcal{N}$ and $\mathcal{N}'$ to be related by a linear transformation and therefore share many of the same qualitative features (e.g. number and stability of equilibria). This phenomenon was termed \emph{linear conjugacy} in \cite{J-S2}.

For completeness, we include the formal definition of linear conjugacy as presented in \cite{J-S2}. We will let $\Phi(\mathbf{x}_0,t)$ denote the flow of (\ref{de}) associated with $\mathcal{N}$ and $\Psi(\mathbf{x}_0,t)$ denote the flow of (\ref{de}) associated with $\mathcal{N}'$.

\begin{definition}
\label{conjugate1}
We will say two chemical reaction networks $\mathcal{N}$ and $\mathcal{N}'$ are \textbf{linearly conjugate} if there exists an injective and surjective linear mapping $\mathbf{h}: \mathbb{R}^n_{>0} \mapsto \mathbb{R}_{>0}^n$ such that $\mathbf{h}(\Phi(\mathbf{x}_0,t))=\Psi(\mathbf{h}(\mathbf{x}_0),t)$ for all $\mathbf{x}_0 \in \mathbb{R}_{>0}^n$.
\end{definition}
\noindent It is known that an injective and surjective linear transformations $\mathbf{h}: \mathbb{R}^n_{>0} \mapsto \mathbb{R}_{>0}^n$ can consist of at most positive scaling and reindexing of coordinates (Lemma 3.1, \cite{J-S2}). Linear conjugacy has been subsequently studied from a computational point of view in \cite{J-S4} and \cite{J-S5}.

Linear conjugacy is a generalization of the concept of \emph{dynamical equivalence} whereby two reaction networks with different topological network structure can generate the same exact set of differential equations (\ref{de}). Two dynamically equivalent networks $\mathcal{N}$ and $\mathcal{N}'$ are said to be alternative \emph{realizations} of the kinetics (\ref{de}), although it is sometimes preferable to say that $\mathcal{N}$ is an alternative realization of $\mathcal{N}'$ or vice-versa. Since the case of two networks being realizations of the same kinetics is encompassed as a special case of linear conjugacy taking the transformation to be the identity, we will focus on linearly conjugate networks.


If a network $\mathcal{N}$ can be shown to be linearly conjugate to a network $\mathcal{N}'$ from a class of networks with desired properties, those properties of $\mathcal{N}'$ are transferred to $\mathcal{N}$. This raises the question of how to find a linearly conjugate network $\mathcal{N}'$ when only the original network $\mathcal{N}$ is given. This was studied in \cite{J-S4} and \cite{J-S5} where the authors built upon the linear programming algorithm introduced in \cite{Sz2}. We can impose that a network $\mathcal{N}'$ be linearly conjugate to a given network $\mathcal{N}$ with the set of linear constraints
\begin{equation}
\label{conjugate}
\mbox{\textbf{(LC)}} \; \; \left\{ \; \; \begin{array}{ll} & \displaystyle{Y \cdot A_b = T^{-1} \cdot M} \\ & \displaystyle{\sum_{i=1}^m [A_b]_{ij} = 0, \; \; \; j=1, \ldots, m} \\ & \displaystyle{0 \leq [A_b]_{ij} \leq 1 / \epsilon, \; \; \; i,j = 1, \ldots, m, \; i \not= j} \\ & \displaystyle{[A_b]_{ii} \leq 0, \; \; \; i = 1, \ldots, m} \\ &\displaystyle{\epsilon  \leq c_j \leq 1/\epsilon, \; \; \; j=1, \ldots, n} \end{array} \right.
\end{equation}
where $0 < \epsilon \ll 1$, and the matrices $M \in \mathbb{R}^{n \times m}$ and $T \in \mathbb{R}^{n \times n}$ are given by:
\begin{eqnarray}
\label{M} M & = & Y \cdot A_k, \mbox{ and} \\
\label{T} T & = & \mbox{diag}\left\{ \mathbf{c} \right\}.
\end{eqnarray}
The kinetics matrix for the network $\mathcal{N}'$ can by constructed from $A_b \in \mathbb{R}^{m \times m}$ and $\mathbf{c} \in \mathbb{R}_{>0}^n$ by the relation
\begin{equation}
\label{newrateconstants}
A_k' = A_b \cdot \mbox{diag} \left\{ \Psi (\mathbf{c}) \right\}.
\end{equation}

Finding a network satisfying (\ref{conjugate}) and then solving (\ref{newrateconstants}) is sufficient to determine a network $\mathcal{N}'$ which is linearly conjugate to $\mathcal{N}$ via the transformation $\mathbf{h}(\mathbf{x}) = T^{-1} \mathbf{x}$. It is important to note that, while the matrix $A_b$ does not exactly correspond to the kinetics matrix $A_k'$ of the network $\mathcal{N}'$, it is a kinetics matrix and has the same structure as $A_k'$. Most importantly, $A_b$ corresponds to a weakly reversible network if and only if the $A_k'$ defined by (\ref{newrateconstants}) does.

It is also worth noting that, without further clarification, the value $m$ may be different in the linear conjugacy constraints (\ref{conjugate}) than that in Definition \ref{deficiency}. The value $m$ in the constraints (\ref{conjugate}) corresponds to the set of \emph{potential} complexes represented in the matrix $Y$. This set must be fixed before the optimization procedure can begin and as such there is no way to guarantee \emph{a priori} that all of the complexes assigned to $Y$ will be present in the network $\mathcal{N}'$. The value may consequently be different than the value $m$ used in Definition \ref{deficiency} to represent the number of complexes in the network $\mathcal{N}'$. (It is common to initialize $Y$ with the reactant and product complexes from the original network $\mathcal{N}$ \cite{Sz2,Sz-H,Sz-H-T,J-S4}; however this choice is by no means necessarily the best \cite{J-S5}.)

We will avoid this notational difficulty by expanding the definition of the network $\mathcal{N}'$ to include all of the complexes initialized in $Y$; that is to say, we will allow complexes to be in $\mathcal{N}'$ even if they are not the reactant or product complex for any reaction in $\mathcal{N}'$. Consequently, the values of $m$ in the linear constraints and the value of $m$ in Definition \ref{deficiency} will coincide for the network $\mathcal{N}'$.

\section{Minimal Deficiency Networks}

We saw in Section \ref{deficiencysection} that many properties of the equilibrium set of a mass-action system (\ref{de}) can be related to the value of the deficiency of the corresponding chemical reaction network. Theorem \ref{deficiencyzerotheorem} and Theorem \ref{conditionstheorem} show that weakly reversible networks with lower deficiency are closer to complex balanced networks (and the regular dynamics guaranteed by Theorem \ref{cbtheorem}) than ones with higher deficiency. Theorem \ref{deficiencyonetheorem} gives conditions under which a network with a nonzero deficiency can still be guaranteed to have a unique positive equilibrium concentration in each compatibility class.

In both of these cases, we find that for weakly reversible networks a \emph{lower} deficiency is more likely to be related to regular dynamics than a \emph{higher} deficiency. In cases where a network is linearly conjugate to multiple weakly reversible networks, therefore, we are likely to be most interested in the network with the minimal deficiency. In this section, we develop a mixed-integer linear programming procedure capable of determining a linearly conjugate weakly reversible network with a minimal deficiency.

\subsection{Parameters of Interest}

We recall that, according to Definition \ref{deficiency}, the deficiency of a chemical reaction network depends on three structural parameters of the reaction network: the number of stoichiometrically distinct complexes $m$, the number of linkage classes $\ell$, and the dimension of the stoichiometric space $s$.

In order to avoid ambiguity with the value $m$ used to represent the number of \emph{potential} complexes in the linear constraints (\ref{conjugate}) and the value in Definition \ref{deficiency}, we allow $\mathcal{N}'$ to contain complexes which do not correspond to the reactant or product complex of a reaction in $\mathcal{N}'$. These unused complexes will appear in the reaction graph as isolated nodes. As such, each unused complex will correspond to a linkage class in itself and this linkage class will be trivially strongly connected. Consequently, including these unused complexes in $\mathcal{N}'$ will not change the value of the deficiency of $\mathcal{N}'$, since the increase to $m$ will be offset by a corresponding increase in $\ell$ in Definition \ref{deficiency}. Including these complexes will also not alter properties related to the reversibility of the network.

In order to minimize the deficiency within a linear programming framework, we need to find a way to quantify the parameters $m$, $\ell$, and $s$. We notice that, since linear conjugacy preserves the dimension of invariant linear spaces, the dimension of the stoichiometric space $s$ is fixed for $\mathcal{N}'$ by the dimension of the largest invariant linear space of the original network $\mathcal{N}$. Since $\mathcal{N}'$ is weakly reversible, the parameter $s$ corresponds to the dimension of the stoichiometric space (see Corollary 1 of \cite{H-F}); however, this correspondence is not necessarily true of $\mathcal{N}$ (see Example 2 of \cite{H-J2}).

We also notice that the value of $m$ is predetermined for the optimization procedure. Consequently, in order to minimize the deficiency it is sufficient to maximize $\ell$. In other words, we need to maximize the number of linkage classes, where each unused complex corresponds trivially to its own linkage class.


\subsection{Counting Linkage Classes}
\label{coolsection}

We need to keep track of the number of linkage classes. In general, this is a difficult task; however, we are aided by the following result.

\begin{theorem}[Theorem 3.1 of \cite{G-H} and Proposition 4.1 of \cite{F3}]
\label{weaklyreversible}
Let $A_k$ be a kinetics matrix and let $\Lambda_i$, $i=1, \ldots, \ell,$ denote the support of the $i^{th}$ linkage class. Then the reaction graph corresponding to $A_k$ is weakly reversible if and only if there is a basis of ker$(A_k)$, $\left\{ \mathbf{b}^{(1)}, \ldots, \mathbf{b}^{(\ell)} \right\}$, such that, for $i=1, \ldots, \ell$,
\begin{equation}
\label{bi}
\mathbf{b}^{(i)} = \left\{ \begin{array}{ll} b^{(i)}_j > 0, \hspace{0.3in} & j \in \Lambda_i \\ b^{(i)}_j = 0, & j \not\in \Lambda_i. \end{array} \right.
\end{equation}
\end{theorem}
\noindent It is easy to see that, for a weakly reversible network, the dimension of ker$(A_k)$ given by Theorem \ref{weaklyreversible} corresponds to the number of linkage classes. Consequently, the parameter $\ell$ here coincides with the earlier usage.

It is typical to assume in Theorem \ref{weaklyreversible} that the kinetics matrix $A_k$ corresponds to a network where every complex appears on either the reactant or product side of at least one reaction. It is easy to extend this to the case where complexes are not used by the reaction network by noting that any unused complexes will contribute an element to ker$(A_k)$ satisfying (\ref{bi}) corresponding to a single positive value in the coordinate corresponding to the unused complex and zeroes elsewhere. In other words, for weakly reversible networks with unused complexes, we can extend the basis of ker$(A_k)$ by considering unused complexes \emph{as their own linkage classes}. Throughout this section, we will allow $\ell$ to correspond to both traditional linkage classes containing several complexes \emph{and} unused complexes.

Theorem \ref{weaklyreversible} implies that for a weakly reversible network, the supports of the basis elements of ker$(A_k)$ represent a complete partition of the set $\left\{ 1, \ldots, m \right\}$, that is to say, we require
\begin{equation}
\label{test}
\begin{array}{c}\Lambda_{k_1} \cap \Lambda_{k_2} = \emptyset, \; \; \mbox{for all } k_1, k_2=1, \ldots, \ell, k_1 \not= k_2 \\
\displaystyle{\bigcup_{k=1}^{\ell} \Lambda_k = \left\{ 1, \ldots, m \right\}}.
\end{array}
\end{equation}

The value of $\ell$, however, is not known; in fact, it is what is to be determined through the procedure. In order to implement this into a computational framework, therefore, we need to determine an upper limit for the number of possible supports $\Lambda_k$. We recall that the deficiency is a nonnegative parameter, so that we have $\delta = m - \ell - s \geq 0$ \cite{F1,H}. Consequently, we have $\ell \leq m - s$ and therefore may use $m-s$ as an upper bound for $\ell$. It is also clear that the deficiency will be zero if and only if this upper bound is attained.

We now introduce the binary variables $\gamma_{ik} \in \left\{0, 1 \right\}$, for $i=1, \ldots, m$, $k=1, \ldots, m-s$, defined according to
\begin{equation}
\label{gamma}
\gamma_{ik} = \left\{ \begin{array}{ll} 1, \; \; \; \; \; & \mbox{if } i \in \Lambda_k \\ 0, & \mbox{if } i \not\in \Lambda_k. \end{array} \right.
\end{equation}
The variables $\gamma_{ik}$ keep track of how the supports of the basis vectors in ker$(A_k)$ according to (\ref{bi}) partition the set $\left\{ 1, \ldots, m \right\}$. For weakly reversible networks, this corresponds to an assignment between the complexes and the linkage classes by Theorem \ref{weaklyreversible}. In other words, $\gamma_{ik} = 1$ if and only if $C_i \in \mathcal{L}_k$.

We also introduce variables $\theta_k \in [0, 1]$, $k=1, \ldots, m-s$, defined according to
\begin{equation}
\label{theta}
\theta_k = \left\{ \begin{array}{ll} 1, \; \; \; \; \; & \mbox{if supp}(\Lambda_k) \not= \emptyset \\ 0, & \mbox{if supp}(\Lambda_k) = \emptyset.\end{array} \right.
\end{equation}
\noindent The variables $\theta_k$ keep track of whether the $k$th partition of $\left\{1, \ldots, m \right\}$ is empty or nonempty. It should be noted that, while we would like the $\theta_k$'s to count the number of non-empty supports, and are therefore interested in only the values $\theta_k=0$ and $\theta_k=1$, it will be possible to relax the integrality of the $\theta_k$'s to vary continuously within the range $[0,1]$. This will be justified in Section \ref{minimizedeficiencysection}.

In order to accommodate the \emph{complete partition} requirements (\ref{test}) as linear constraints, and to accommodate (\ref{theta}), we impose
\begin{equation}
\label{partition}
\mbox{\textbf{(CP)}} \; \; \left\{ \; \; \begin{array}{ll} & \displaystyle{\sum_{k=1}^{m-s} \gamma_{ik} = 1, \; \; \; i=1, \ldots, m} \\ & \displaystyle{\sum_{i=1}^m \gamma_{ik} - \epsilon \theta_k \geq 0, \; \; \; k = 1, \ldots, m-s} \\ & \displaystyle{-\sum_{k=1}^m \gamma_{ik} + \frac{1}{\epsilon} \theta_k \geq 0, \; \; \; k = 1, \ldots, m-s} \\ & \gamma_{ik} \in \left\{ 0, 1 \right\}, \; \; \; i = 1, \ldots, m, k = 1, \ldots, m-s \\ & \theta_k \in \left[ 0, 1 \right], \; \; \; k = 1, \ldots, m-s. \end{array} \right.
\end{equation}
\noindent where $0 < \epsilon \ll 1$ is sufficiently small and can be chosen to be the same as in (\ref{conjugate}).

The first constraint set guarantees that each complex appears in exactly one partition. The second two constraint sets of (\ref{partition}) correspond to the constraint
\[0 \leq \epsilon \theta_k \leq \sum_{i=1}^m \gamma_{ik} \leq \frac{1}{\epsilon} \theta_k\]
which keeps track of whether the $k$th partition is empty or nonempty. If no element is in the $k$th partition, then the sum is zero, which forces $\theta_k$ to be zero (first inequality). If there is an element in the $k$th partition, then the sum is nonzero, which forces $\theta_k$ to be nonzero (second inequality). An argument in Section \ref{minimizedeficiencysection} will allow us to conclude that any nonzero $\theta_k$ must be one, so that this fulfills the requirements of (\ref{theta}).

\subsection{Constructing the Kernel}

We still need to guarantee that the sets $\Lambda_k$ considered in Section \ref{coolsection} correspond to the supports of vectors in ker$(A_b)$ which satisfy (\ref{bi}). In other words, we need to restrict ourselves to sets $\Lambda_i$, $i=1, \ldots, \ell$, where there exists a $\mathbf{b}^{(i)}$ satisfing (\ref{bi}) and
\begin{equation}
\label{21}
A_b \cdot \mathbf{b}^{(i)} = \mathbf{0}, \; \; \; i = 1, \ldots, \ell.
\end{equation}

We follow the technique outlined in the paper \cite{J-S4} for determining weakly reversible networks. We define a matrix $\Phi \in \mathbb{R}^{m \times m}$ with entries
\begin{equation}
\label{construct}
\Phi_{ij} = [A_b]_{ij} \cdot \mathbf{b}_j
\end{equation}
where $\mathbf{b} = \sum_{k=1}^\ell \mathbf{b}^{(k)}$ and $\left\{ \mathbf{b}^{(1)}, \ldots, \mathbf{b}^{(\ell)} \right\}$ is any set of vectors satisfying (\ref{bi}) and forming a basis of ker$(A_b)$. We can see that the system of non-linear equations (\ref{21}) is satisfied if and only if the system of linear equations
\begin{equation}
\label{kernelsum}
\sum_{j \in \Lambda_k} \Phi_{ij} = 0, \; \; \; \mbox{for all } i = 1, \ldots, m, k=1, \ldots, \ell,
\end{equation}
is satisfied.

We know that $\Phi$ is a kinetics matrix since $A_b$ is and (\ref{construct}) preserves this property. That is to say, we have
\begin{equation}
\label{kineticproperty}
\displaystyle{\sum_{j=1}^m \Phi_{ji} = 0, \; \; \; \mbox{for all } i = 1, \ldots, m}
\end{equation}
and
\[\displaystyle{\Phi_{ij} \geq 0, \; \; \; \mbox{for all } i, j =1, \ldots, m, i \not= j.}\]
Consequently, we can solve the diagonal elements of $\Phi$ in (\ref{kineticproperty}) and substitute them into (\ref{kernelsum}) to get the simplified constraint set
\begin{equation}
\label{35}
\mathop{\sum_{j \in \Lambda_k}}_{j \not= i} \Phi_{ij} = \mathop{\sum_{j=1}^m}_{j \not= i} \Phi_{ji}, \; \; \; \mbox{for } i = 1, \ldots, m, k = 1, \ldots, \ell.
\end{equation}

We need to derive linear constraints capable of constructing a matrix $\Phi$ according to (\ref{construct}) which satisfies (\ref{35}) and has the same structure as $A_b$. This is made more challenging than the case considered in \cite{J-S4} by the requirement that the kernel vector $\mathbf{b}$ in (\ref{construct}) decompose according to the partitions $\Lambda_k$, $k=1, \ldots, \ell$.

We notice first of all that we do not know how many partitions $\Lambda_k$ we will need, so we will take the upper bound $m-s$ on the number of partitions (see Section \ref{coolsection}). We also notice that the construction (\ref{construct}) requires that $\Phi_{ij} = 0$ for any $i$ and $j$ such that $i \in \Lambda_k$ and $j \not\in \Lambda_k$ for some $k = 1, \ldots, m-s$. That is to say, we do not permit reactions to occur between complexes in different linkage classes. This can be accommodated by noticing that, if two indices $i$ and $j$ are on the same support $\Lambda_k$, then we will have $\gamma_{ik}-\gamma_{jk} = 0$ \emph{for all} $k = 1, \ldots, m-s$, where as if $i$ and $j$ are not on the same support $\Lambda_k$, we will have $\gamma_{ik}-\gamma_{jk} \in \left\{ -1, 0, 1 \right\}$ and attain the value $-1$ for at least one $k$. Consequently, we can accommodate this requirement with the linear constraints 
\begin{equation}
\label{324}
\begin{array}{c} \displaystyle{\Phi_{ij} \leq 1/\epsilon (\gamma_{ik}-\gamma_{jk}+1),} \vspace{0.05in} \\ \mbox{for } i, j = 1, \ldots, m, i \not= j, k=1, \ldots, m-s. \end{array}
\end{equation}

We can accommodate (\ref{35}), (\ref{324}), and the requirement that $\Phi$ have the same structure as $A_b$, with the linear constraint set
\begin{equation}
\label{kernel}
\mbox{\textbf{(Ker)}} \; \; \left\{ \; \; \begin{array}{ll} & \displaystyle{\mathop{\sum_{l=1}^m}_{l \not= i} \Phi_{il} = \mathop{\sum_{l=1}^m}_{l \not= i} \Phi_{li}} \\ & \displaystyle{\Phi_{ij} \leq \frac{1}{\epsilon} (\gamma_{ik}-\gamma_{jk}+1)} \vspace{0.05in} \\ & \Phi_{ij} \geq \displaystyle{\epsilon [A_b]_{ij}} \vspace{0.05in} \\ & \Phi_{ij} \leq \displaystyle{\frac{1}{\epsilon}[A_b]_{ij}} \vspace{0.05in} \\ & i,j = 1, \ldots, m, i \not= j, \; \; k = 1, \ldots, m-s. \end{array} \right.
\end{equation}
Notice that (\ref{35}) can be generalized to the first constraint in (\ref{kernel}) because of the imposition that $\Phi_{il}=0$ for any $l \not\in \Lambda_k$ in the left sum. This is guaranteed by (\ref{324}).

\subsection{Uniqueness of Solution}
\label{uniquenesssection}

The constraint sets (\ref{conjugate}), (\ref{partition}), and (\ref{kernel}), form the basis of a mixed-integer linear programming (MILP) problem. This is due to the non-continuous binary variables $\gamma_{ik}$, $i=1, \ldots, m$, $j=1, \ldots, m-s,$ which keep track of how the complexes are assigned to the partitions.

MILP problems are known to NP-hard and are generally approached by a branch-and-bound method (for an accessible introduction to branch-and-bound methodology, see \cite{Fl,N-W}). One well-known complicating factor for branch-and-bound methods is non-uniqueness of the integer portion of the problem. There is nothing in the constraint sets (\ref{partition}) which guarantee a unique assignment of the complexes to partitions. That is to say, a linkage class could be assigned to the first partition as easily as the second or third. Consequently, we would like to introduce further constraints which guarantee a unique partitioning structure for the optimal solution.

An intuitive way to structure the partitioning variables $\gamma_{ik}$ is to always assign the first complex to the first partition, and then to assign each subsequent free complex to the next available partition. That is to say, if the second complex is not in the same partition as the first complex, it will be assigned to the second partition. If it is in the same partition as the first complexes, but the third is not, then the third complex will be assigned to the second partition.

This \emph{unique} partition structure can be guaranteed by the linear constraint set
\begin{equation}
\label{unique}
\mbox{\textbf{(U)}} \; \; \left\{ \; \; \begin{array}{ll} & \displaystyle{\sum_{j = 1}^{i-1} \gamma_{jk} \geq \sum_{l = k+1}^{m-s} \gamma_{il},} \\ & i = 1, \ldots, m, k=1, \ldots, m-s, k \leq i.\end{array} \right.
\end{equation}

This constraint set guarantees that, in ascending order, each complex which is assigned to a new partition is assigned to the unused partition with the lowest index. It does so in the following way:
\begin{enumerate}
\item
For the case $i=k=1$, the left sum in (\ref{unique}) is empty and therefore returns a value of zero. This forces $\gamma_{1l}=0$ for $l=2, \ldots, m-s$, which forces $\gamma_{11}=1$ by the first condition of (\ref{partition}). In other words, the first complex is always assigned to the first partition.
\item
The left sum in (\ref{unique}) is zero for $i=k=2$, which forces $\gamma_{2j}=0$ for $j=3, \ldots, m-s$. This forces the second complex to be assigned to either the first partition or the second partition by (\ref{partition}), depending on whether it is grouped with the first complex or not.
\item
Following an induction on the complexes as $i=3, \ldots, m$, we can see that the constraint (\ref{unique}) corresponding to the $i$th complex and the first unused partition will always guarantee this complex may not be assigned to a partition of a higher index than the first unused partition. Consequently, if it not grouped with an earlier complex, it must be assigned the index of the next available partition.
\item
It is easy to see that the conditions (\ref{unique}) are satisfied for all of the entries not corresponding to leading ones, so that we are done.
\end{enumerate}


\subsection{Minimizing the Deficiency}
\label{minimizedeficiencysection}

We know that minimizing the deficiency $\delta$ is equivalent to maximizing the number of linkage classes $\ell$, where we include the unused complexes as linkage classes unto themselves.

In Section \ref{coolsection} we introduced variables $\theta_k \in \left[ 0, 1 \right], k=1, \ldots, m-s,$ according to (\ref{theta}) to keep track of the number of partitions $\Lambda_k$ of $\left\{ 1, \ldots, m \right\}$ which corresponded to vectors in ker($A_b$) satisfying $(\ref{bi})$. We have no guarantee, however, that these sets correspond to a complete basis of ker$(A_b)$ since any two vectors in ker($A_b$) satisfying (\ref{bi}) can be added to one another to produce another vector in ker$(A_b)$ also satisfying (\ref{bi}). (This corresponds to multiple linkage classes being placed on the same support $\Lambda_k$.) We also have no guarantee that the $\theta_k$'s will properly enumerate the \emph{number} of basis elements of ker$(A_b)$ since any $\theta_k$ corresponding to a nonempty support may attain a value anywhere between zero and one, and not just the value of one.

We need to guarantee that the sum of the $\theta_k$'s corresponds to the maximal partition of $\left\{ 1, \ldots, m \right\}$ and that each $\theta_k$ attains a value of one when it corresponds to a nonempty support. Consider the objective function corresponding to maximizing the number of linkage classes, which corresponds to minimizing the deficiency. We have
\begin{equation}
\label{minimizedeficiency}
\mbox{\textbf{(Min Def)}} \; \; \left\{ \; \; \; \; \; \; \; \; \mbox{minimize} \; \; \; \; \; \sum_{k=1}^{m-1} -\theta_k \right.
\end{equation}
over the constraint sets (\ref{conjugate}), (\ref{partition}), and (\ref{kernel}). We notice that this objective function guarantees that, if a $\theta_k$ is in between zero and one, it will attain the value one and that a maximal partition of $\left\{ 1, \ldots, m \right\}$ will be chosen, because in both cases such a situation is more optimal than the alternative. In other words, the $\theta_k$'s properly count the number of linkage classes when (\ref{minimizedeficiency}) is imposed.

If the off-diagonal elements are solved for in (\ref{conjugate}), the algorithm presented here for finding a weakly reversible network with minimal deficiency which is linearly conjugate to a given network contains $m(2m-1)+n-s$ continuous decision variables and $m(m-s)$ binary decisions variables.

It should be noted that a linearly conjugate network with the maximal deficiency cannot be obtained in the same fashion as outlined here. This is due to having to maximize the sum in (\ref{minimizedeficiency}) in order to make the correspondence between the sum of the $\theta_k$'s and the value of $\ell$. Without this step, we could not rule out linkage classes clumping together onto the same support in ker$(A_k)$.

\section{Examples}

We now introduce a few examples which illustrate the methods presented in this paper. All computations are performed on the primary author's personal-use Acer laptop (AMD Athlon II Neo K125 Processor 1.70 GHz, 4 GB RAM).\\

\begin{figure}[h]
\begin{center}
\includegraphics[width=7cm]{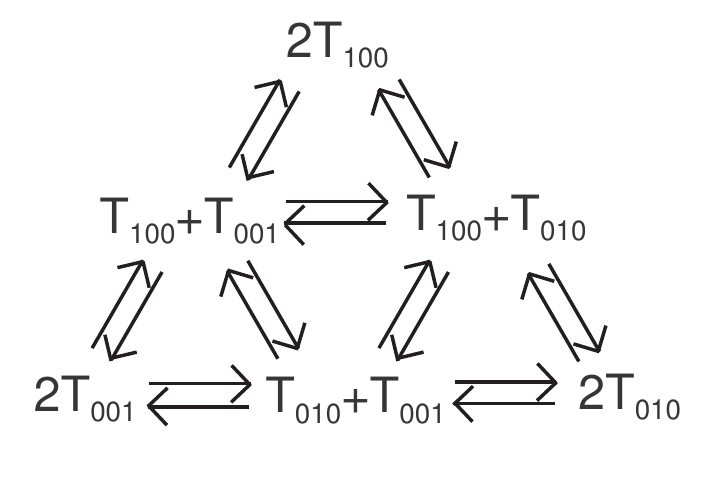}
\end{center}
\vspace{-0.3in}
\caption{Network of singularly bound enzymes with three binding sites. Interactions between the enzymes allow transfer of substrates from one binding site to another.}
\label{figure4}
\end{figure}

\noindent \textbf{Example 1:} We propose the following mechanism. Consider a substrate which can bind to an enzyme $T$ at one of three binding sites and let $T_{100}$, $T_{010}$, and $T_{001}$ denote the enzyme with the substrate bound at the first, second, and third binding site, respectively. Suppose that binary collisions between the substrates can cause a spontaneous shift in the substrate from one binding site to another, but no transfer of substrate from one bound enzyme to another. This web of interactions can be visualized by the network in Figure \ref{figure4}.

This network is weakly reversible and has a deficiency of three ($m=6, \ell=1, s=2$, $\delta = m - \ell - s = 3$). Due to the high deficiency, we may not apply Theorem \ref{deficiencyzerotheorem} or Theorem \ref{deficiencyonetheorem}; however, since the network is weakly reversible, we may apply Theorem \ref{conditionstheorem}. This tells us that there are three algebraically independent conditions on the rate constants which must be satisfied in order for the system to be complex balanced and therefore fall within the scope of Theorem \ref{cbtheorem}. If we set $C_1=2T_{100}$, $C_2=2T_{010}$, $C_3=2T_{001}$, $C_4=T_{100}+T_{010}$, $C_5=T_{100}+T_{001}$, and $C_6=T_{010}+T_{001}$, the required conditions are $K_1 K_2 = K_4^2$, $K_1 K_3 = K_5^2$, and $K_2 K_3 = K_6^2$, where $K_i$ is the absolute value of the $(i \times i)^{th}$ minor of $A_k$. (For details, see \cite{C-D-S-S}.)

We might wonder what additional behaviour is permitted by the reaction network given in Figure \ref{figure4}. A general analysis can be conducted by the CRN toolbox made freely available online \cite{J-E-K-F}. This is a powerful computational toolbox capable of determining whether a chemical reaction network has the capacity for zero eigenvalues, injectivity, multistability, and concordance \cite{Fe,Fe2,E-F,Sh-F2}. The toolbox's \emph{Higher Deficiency Report} reveals that the mass-action system corresponding to the network has the capacity for multiple positive steady states for the rate constant choices
\begin{equation}
\label{rateconstant1}
\begin{array}{lll} k(1,4) = 7.389056 & k(4,1) = 2.7182818 & k(5,4) = 4.3002585 \\ k(1,5) = 7.389056 & k(4,2) = 2.7182818 & k(5,6) = 4.3002585 \\ k(2,4) = 1 & k(4,5) = 55.125832 & k(6,2) = 32.08195 \\ k(2,6) = 1 & k(4,6) = 29.019118 & k(6,3) = 1.5819767 \\ k(3,5) = 2.5026503 & k(5,1) = 45.90757 & k(6,4) = 1.5819767 \\ k(3,6) = 2.5026503 & k(5,3) = 4.3002585 & k(6,5) = 1.5819767.\end{array}
\end{equation}

For rate constant choices lying away from these sets, the dynamics is significantly more difficult to analyze. We consider the deficiency-reducing methodology introduced in this paper on two sets of rate constant choices (in addition to (\ref{rateconstant1})). The two additional sets of rate constant choices are
\begin{equation}
\label{rateconstant2}
k(i,j)=j, \; \; \; i,j=1, \ldots 6, \; \; \; (i,j) \in \mathcal{R}
\end{equation}
and
\begin{equation}
\label{rateconstant3}
k(i,j)=i, \; \; \; i,j=1, \ldots, 6, \; \; \; (i,j) \in \mathcal{R}. 
\end{equation}
In other words, we consider the rate constant choices where all the rates going into a particular complex are the same, and the rate constant choices where all the rates going out of a particular complex are the same. For simplicity, we take the differences between complexes to scale according to the index of the complex itself.

We optimize (\ref{minimizedeficiency}) over the contraint sets (\ref{conjugate}), (\ref{partition}), (\ref{kernel}), and (\ref{unique}), in GLPK. The algorithm runs in under a second for all rate constant sets. The rate constants choice (\ref{rateconstant1}) produces a network with deficiency three, and consequently we conclude that there is no network which is linearly conjugate to the network in Figure \ref{figure4} with a lower deficiency.

\begin{figure}[h]
\begin{center}
\includegraphics[width=11cm]{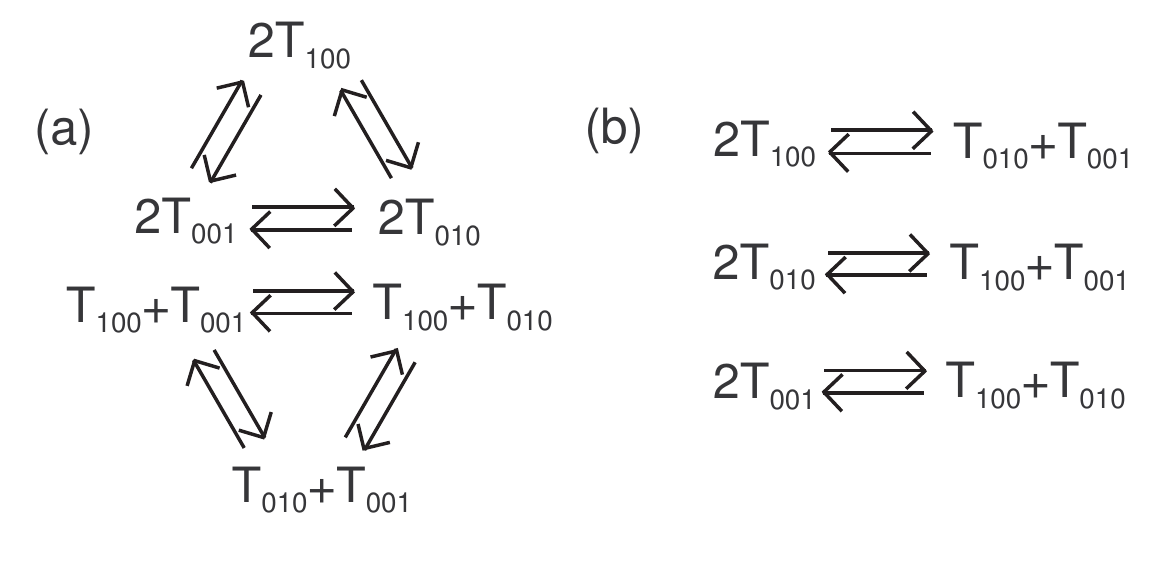}
\end{center}
\vspace{-0.3in}
\caption{Weakly reversible networks which are dynamically equivalent to the one contained in Figure \ref{figure4}. The network in (a) is dynamically equivalent for the rate constant choices (\ref{rateconstant2}) and has a deficiency of two. The network in (b) is dynamically equivalent for the rate constant choices (\ref{rateconstant3}) and has a deficiency of one.}
\label{figure5}
\end{figure}

The results of the optimization for (\ref{rateconstant2}) and (\ref{rateconstant3}) are contained in Figure \ref{figure5}. In both cases, the conjugacy constants are $c_1=c_2=c_3=1$. The network corresponding to the rate constant set (\ref{rateconstant2}) is deficiency two while the network corresponding to the rate constant set (\ref{rateconstant3}) is deficiency one. Neither network is amenable to application of Theorem \ref{deficiencyzerotheorem} or Theorem \ref{deficiencyonetheorem}. The CRN toolbox can, however, by used to verify (by the algorithm presented in \cite{Fe4}) that the network contained in Figure \ref{figure5}(b) admits at most one equilibrium concentration in each positive compatibility class. (Since weakly reversible networks contain at least one equilibrium concentration in each compatibility class, this is sufficient to guarantee that the network has exactly one equilibrium concentration in each compatibility class \cite{D-F-J-N}.)

The rate constant values corresponding to the two networks in Figure \ref{figure5} are
\[\begin{array}{lllllll} (a) & k(1,2) = 2 & k(3,1) = 2.5 & k(5,4) = 2 \\ & k(1,3) = 2.5 & k(3,2) = 3 & k(5,6) = 3 \\ & k(2,1) = 2 & k(4,5) = 4 & k(6,4) = 3 \\ & k(2,3) = 3 & k(4,6) = 7 & k(6,5) = 6, \\ & & & \\ (b) & k(1,6) = 1 & k(3,4) = 3 & k(5,2) = 5 \\ & k(2,5) = 2 & k(4,3) = 4 & k(6,1) = 6. \end{array}\]

\noindent \textbf{Example 2:} Consider the set of differential equations given by
\begin{equation}
\label{example2}
\begin{split}
\frac{dx_1}{dt} & = 1-x_1^2-x_1+x_2x_3 \\ \frac{dx_2}{dt} & = 2x_1-2x_2x_3-2x_2^2+2x_3^2 \\ \frac{dx_3}{dt} & = x_1-x_2x_3+x_2^2-x_3^2.
\end{split}
\end{equation}
We can produce a chemical reaction network which generates (\ref{example2}) using the algorithm presented in \cite{H-T} and reproduced in \cite{Sz-H}. This gives a network involving the complexes
\[\begin{split} & C_1 = \emptyset, C_2 = X_1, C_3 = 2X_1, C_4 = X_2+X_3,\\ & C_5 = X_1+X_2+X_3, C_6 = X_1+X_2, C_7 = X_3, \\ & C_8 = 2X_2, C_9 = X_2, C_{10} = 2X_3, C_{11} = X_2 + 2 X_3,\\ & C_{12} = X_1+X_3, C_{13} = 2X_2+X_3.\end{split}\]

The corresponding network is not weakly reversible and has a deficiency of eight. It is therefore not amenable to Theorem \ref{conditionstheorem} or Theorem \ref{deficiencyonetheorem}. In order to apply these results, therefore, we would like to find a weakly reversible network which is linearly conjugate to this one and has a minimal deficiency. We construct the matrices
\[Y = \left[ \begin{array}{ccccccccccccc} 0 & 1 & 2 & 0 & 1 & 1 & 0 & 0 & 0 & 0 & 0 & 1 & 0 \\ 0 & 0 & 0 & 1 & 1 & 1 & 0 & 2 & 1 & 0 & 1 & 0 & 2 \\ 0 & 0 & 0 & 1 & 1 & 0 & 1 & 0 & 0 & 2 & 2 & 1 & 1 \end{array} \right] \]
and
\[M = \left[ \begin{array}{ccccccccccccc} 1 & -1 & -1 & 1 & 0 & 0 & 0 & 0 & 0 & 0 & 0 & 0 & 0 \\ 0 & 2 & 0 & -2 & 0 & 0 & 0 & -2 & 0 & 2 & 0 & 0 & 0 \\ 0 & 1 & 0 & -1 & 0 & 0 & 0 & 1 & 0 & -1 & 0 & 0 & 0 \end{array} \right]\]
and optimize (\ref{minimizedeficiency}) over the contraint sets (\ref{conjugate}), (\ref{partition}), (\ref{kernel}), and (\ref{unique}). The upper limit for the number of partitions, $m-s$, can be easily determined to be $10$.

\begin{figure}[h]
\begin{center}
\includegraphics[width=11cm]{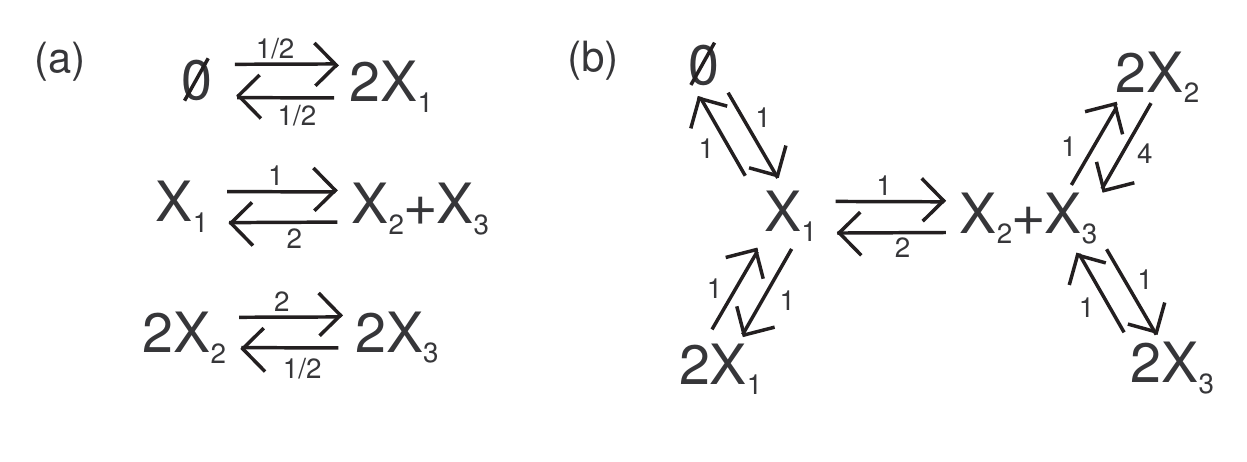}
\end{center}
\vspace{-0.3in}
\caption{Two weakly reversible network structures which are linearly conjugate to a network generating kinetics (\ref{example2}). The conjugacy constants are $c_1=c_3=1$ and $c_2=2$. The network in part (a) has a deficiency of zero while the network in part (b) has a deficiency of two. Isolated complexes have been excluded from the reaction graph.}
\label{figure2}
\end{figure}

A quick computation in GLPK produces the network given in Figure \ref{figure2}(a). We can readily see that this is a weakly reversible zero deficiency system and therefore falls within the scope of the networks considered by Theorem \ref{deficiencyzerotheorem} and Theorem \ref{cbtheorem}. Consequently, we know that (\ref{example2}) has exactly one positive equilibrium concentration and that this equilibrium concentration is locally asymptotically stable.

It is interesting to note that the network contained in Figure \ref{figure2}(a) is not the only weakly reversible network capable of producing (\ref{example2}). If we do not insist on maximizing the number of linkage classes, other networks can be selected with a sub-optimal deficiency value. The network in Figure \ref{figure2}(b) is also linearly conjugate to a network generating (\ref{example2}) with the same conjugacy constants as that of Figure \ref{figure2}(a). Despite being weakly reversible, however, this network has a deficiency of two and is not amenable to either Theorem \ref{deficiencyzerotheorem} or Theorem \ref{deficiencyonetheorem}.

\section{Conclusions}

In this paper, we have presented a computational method for determining weakly reversible networks with minimal deficiency which are linearly conjugate to a given chemical reaction network or kinetic system.

It was shown that, for this purpose, it is sufficient to maximize the number of linkage classes where the linkage classes are defined to include isolated complexes as well as traditional linkage classes with multiple complexes. The proposed algorithm is based on mixed integer linear programming where the binary variables are used to keep track of the complexes' positions in different linkage classes, and the continuous variables keep track of the reaction rate coefficients, the conjugacy coefficients, the structure of the network, and the emptiness/non-emptiness of complex partitions. An additional linear constraint on the binary variables ensures the uniqueness of the complex partitioning, which dramatically improves the computational efficiency of the algorithm. We then applied the algorithm to several examples and were able to use the Deficiency Zero Theorem (Theorem \ref{deficiencyzerotheorem}) and Deficiency One Theorem (Theorem \ref{deficiencyonetheorem}) to determine properties of the equilibrium set of the original reaction network.

Future work in this area includes:
\begin{enumerate}
\item
In order to apply the algorithm outlined in this paper, it is necessary that the rate constants of the original reaction network be specified. Consequently, we may be overlooking networks with a lower deficiency which are linearly conjugate to a given network for certain rate constants values but not for others. Example 1 provides a good example of a mechanism which permits different optimal deficiencies for different rate constant choices. If we were, for instance, only given the rate constant values specified by (\ref{rateconstant1}) or (\ref{rateconstant2}), we would not realize that the mechanism permitted a deficiency one conjugacy for different rate constant values. (Research on these `structurally-fixed' networks was initiated in \cite{J-S5} but can thus far only solve for dynamically equivalent networks.)
\item
Many deficiency results do not require the networks in consideration to be weak reversibility. For instance, the Deficiency One Theorem (Theorem \ref{deficiencyonetheorem}) and the algorithm presented in \cite{Fe4} only require the networks to have a single terminal strongly linked component within each linkage class. It would be useful, therefore, to adapt the algorithm to networks which are not necessarily weakly reversible but satisfy some alternative structural requirements.
\end{enumerate}

\noindent \textbf{Acknowledgements:} M. Johnston and D. Siegel acknowledge the support of D. Siegel's Natural Sciences and Engineering Research Council of Canada Discovery Grant.  G. Szederk\'{e}nyi acknowledges the support of the Hungarian National Research Fund through grant no. OTKA K-83440 as well as the support of the Marginal Research Association of Budapest.

\bibliographystyle{plain}
\bibliography{myrefs}

\end{document}